\documentclass{amsart}
\usepackage{amssymb}
\usepackage{amsmath}
\usepackage[mathscr]{eucal}
\usepackage{epsfig}
\usepackage{graphicx}
\usepackage{epic}
\usepackage{tikz}
\usepackage{xcolor}

\newtheorem{mytheorem}{\bf Theorem}[section]

\newtheorem{remark}[mytheorem]{Remark}

\newcommand{\cal}{\mathcal }

\newcommand{\e}{\epsilon}

\begin{document}
\title[Porous Media  Fracture Modeling]{ A porous media  fracture model based on homogenization theory }

\author{  J. Galvis }
\address{Departamento de Matem\'aticas, Universidad Nacional de Colombia, Carrera 45 No. 26-85, Edificio Uriel Gutierr\'ez, Bogot\'a D.C., Colombia, \tt jcgalvisa\@unal.edu.co.}
\author{ H. M. Versieux }
\address{Departamento de Matem\'atica, Universidade Federal de Minas Gerais, Belo Horizonte, MG, Brasil, \tt  henrique@im.ufrj.br.}


\date{}

\begin{abstract}
  A novel regularized fracture model for crack propagation in porous media is proposed. Our model is obtained through homogenization theory and formal asymptotic expansions. We start with a regularized quasi-static fracture model posed in a periodically perforated domain obtained by periodic extension of a re-scaled unit cell with a hole. This setup allows us to write two separated minimality conditions for the primary (displacement) and secondary variables plus a balance of energy relation. Then we apply the usual asymptotic expansion matching to deduce limit relations when the re-scaling parameter of the unit cells vanishes. By introducing cell problems solutions and a homogenized tensor we can recast the obtained relations into a novel model for crack propagation in porous media. The proposed model can be interpreted as a regularized quasi-static fracture model for porous media. This model yields two separated (homogenized) minimality conditions for the primary and secondary variables and a balance of homogenized energy relation. 
\end{abstract}

\maketitle


\section{Introduction}

 The   study  of crack propagation models has received considerable attention from the scientific community   recently. For instance, a considerable effort has been invested in the  Francfort and Marigo \cite{MR1633984} model for quasi-static brittle fracture evolution; see for example \cite{MR1897378,MR1978582,MR1988896,MR2186036} and references therein. These studies have focused on a fracture propagating on a given material subject to  loads or  prescribed boundary conditions displacements.   To the best of our knowledge, the modeling of  fracture propagation models on a porous a media  subject to loads or  prescribed boundary conditions displacements have not received considerable attention. For instance, one of the few references found by the authors in this subject is the work   \cite{MR2588959}, where the homogenization of a sharp interface fracture model applied to perforated domains  is studied   in the context of  $\Gamma$-convergence.  The understanding of fracture propagation in porous media  has   important  applications in the oil and gas industry.

In this work we study through formal asymptotic expansions the homogenization of a regularized fracture model.  This regularized quasi static model  was obtained in \cite{MR3460102} in the limit of vanishing viscosity and inertia terms of a dynamic fracture model proposed by  Bourdin, Larsen and Richardson \cite{BoLaRi}. This  model  is   similar to the regularized  version of the Francfort and Marigo model proposed in \cite{MR2390547}. The  difference between these two models concern the minimality condition imposed on the problem, while the regularized Francfort and Marigo model uses a global minimality condition ( see condition $ (\tilde{c}) $ in section \ref{sec:models} below) the regularized model from  \cite{MR3460102}  has a separate minimality condition (see condition $ ({c}) $ in section \ref{sec:models} below). We observe  that in general the numerical approximations of the regularized Francfort and Marigo model correspond to approximations of the model satisfying the separate minimality condition; see \cite{MR3460102}  for more details.

To obtain our porous media fracture model, we consider the regularized fracture model from \cite{MR3460102}  applied to  a perforated  periodic material. We  assume that the periodic material is represented by a domain $\Omega_\epsilon$ given by the 
periodic extension of a scaled  set $\epsilon Y$ over  $\mathbb{R}^d$ contained in a 
open fixed domain $\Omega \subset \mathbb{R}^d$. Here  $Y$ is a reference domain  
containing holes and  $\epsilon>0$ is the  scaling parameter.    Our model is obtained in 
the limit $\epsilon \rightarrow 0$ of the solutions of the brittle  fracture model applied in 
the domain $\Omega_\epsilon$.   We employ asymptotic expansions and formal calculations to derive our model. In particular we work out asymptotic expansions related to 
both separated minimality conditions and combine them with the asymptotic expansions of the energy balance of the model in the auxiliary time variable that the model introduces. The advantage of working with asymptotic expansions is that it allows us to obtain a complete characterization of the homogenized problem in terms of the  local problems associated with the micro-structure geometry. 
Similar strategy has been previously used to derive 
other porous media models  such as the D'Arcy law; see for instance 
 \cite{MR1434316} and references therein.   

\begin{remark} In our work we have two different parameters related to two different physical scales: $(i)$ the $\e$ parameter related to the size of the porous  in the media, $(ii)$ the $\gamma$ regularization parameter  related to the width of the fracture.  Here, we study the case when $\e \rightarrow 0$ and $\gamma$ remains fixed. This framework corresponds to the physical case where the porous size are much smaller than the width of the fracture. 
\end{remark}

The rest of the paper is organized as follows. 
At the end of this section the reader can find some important notation we will use in the rest of the paper. 
In Section \ref{sec:models} we present the regularized 
brittle fracture propagation model we consider  a baseline for our study. In Section \ref{sec:expansions} we introduce the perforated domain set up and the derivations led us to from the asymptotic expansions we consider. In Section
\ref{sec:homogenization} we present our homogenization analysis of the brittle fracture model and in Section \ref{novel} we summarize the obtained novel homogenized model. We round up the paper with some conclusions and perspetive. 

Finally, we  introduce some notation used in our work. We use $\langle \cdot, \cdot \rangle$ to denote the $L^2(\Omega)$ (or  $L^2(\Omega)^n$) inner product or the pairing between $H^1_0(\Omega)$ and $H^{-1}(\Omega)$. We denote the $L^2(\Omega)$ norm  by $\|\cdot\|$, and  the $W^{k,p}(\Omega)$ norm (seminorm)  by  $\|\cdot \|_{k,p}$ ($|\cdot|_{k,p}$). Also, in the case $p=2$ we use the notation $\|\cdot \|_{k}$ ($|\cdot|_{k}$) for the $H^k(\Omega)$ norms (seminorms). The spaces $L^p(0,T;H^1(\Omega))$ are denoted $L^p(0,T;H^1)$ with similar notation 
for other cases of spaces of functions of time and space variable.  The subindex $\epsilon$ represents a sequence (or a subsequence) $\{\epsilon_j\}$ converging to zero. We denote  a generic constant independent of $\epsilon$ by $c$.

\section{Regularized brittle fracture propagation model }\label{sec:models}

We now present the regularized quasi-static fracture model considered in this work. It was obtained in \cite{MR3460102} in the limit of vanishing viscosity and inertia terms of a dynamic fracture model proposed by  Bourdin, Larsen and Richardson \cite{BoLaRi}. For simplicity of the presentation and in order to fix ideas we consider the antiplane case; other cases and models can be considered as well.
In particular we assume that  the displacement of a given material is represented by a function  $u(x):\Omega \rightarrow \mathbb{R}$,  satisfying  a prescribed boundary condition $u(x)|_{\partial \Omega} = g(x)$. Here $\Omega \subset \mathbb{R}^2$ is a bounded open set with Lipschitz continuous boundary.

  First, we  define  the following functionals
\begin{eqnarray}\label{eng_func_E_H}
{\cal E}(u,v)&=& \frac{1}{2} \int_{\Omega} (v^2+ \eta_\gamma) |\nabla u|^2dx,\\\nonumber  \quad {\cal H}(v)&=& \int_{\Omega} \frac{1}{4\gamma} (1-v)^2 + \gamma |\nabla v|^2 dx
\end{eqnarray}
and the total energy as
\begin{equation}\label{ener-func_g}
E_\gamma(u,v)= {\cal E}(u,v) + {\cal H}(v).
\end{equation}
Here $v$ is a  function satisfying $0 \leq v\leq 1$, and  $0 < \eta_\gamma  \ll \gamma$ are the  parameters of  the Ambrosio-Tortorelli approximation $E_\gamma(u,u)$ of the   functional $E(u)$ defined in \cite{MR1633984}, in the sense that,  $E_\gamma$, $\Gamma$-converges to $E$ when $\gamma \rightarrow 0$; see \cite{MR2106765,MR713808,MR1968440}.

The model is given by functions $s \mapsto  u(s)=u(\cdot ,s),\;v(s)=v(\cdot, s)$  such that for every $s\in [0,1]$ we have
\begin{description}
\item[$({a})$]  ${u}(s)=g(s)$ on $\partial \Omega$ and  $0 \leq {v}(s) \leq 1$ for all $s\in [0,1]$ ;
\item[$({b})$]  for all $0 \leq s' \leq s \leq 1$  we have $ {v}(s)\leq {v}(s')$;
\item[$({c})$]  for $s \in [0,1]$ we have the following minimality condition with respect to the secondary variable $v$,
\begin{equation}\label{c.1}
{\cal E}({u}(s),{v}(s)) + {\cal H}({v}(s))= \inf_{0 \leq z \leq {v}(s)}  {\cal E}({u}(s),z) + {\cal H}(z),
\end{equation}
and the following minimality condition with respect to $u$,
\begin{equation}\label{c.2}
{\cal E}({u}(s),{v}(s)) + {\cal H}({v}(s))= \inf_{\phi-g(s) \in H^1_0(\Omega)}  {\cal E}(\phi,{v}(s)) + {\cal H}({v}(s)).
\end{equation}
\item[$({d})$] the function ${\cal E }({u}(s), {v}(s)) + {\cal H}({v}(s))$ is absolutely continuous for  $s \in [0,1]$, and we have the following balance of energy relation
\begin{equation}\label{d}
{\cal E }({u}(s), {v}(s)) + {\cal H}({v}(s)) = {\cal E }({u}(0), {v}(0)) + {\cal H}({v}(0)) +  \int_0^s \langle  (\eta_\gamma + {v}^2) \nabla {u}, \nabla g_s \rangle d\tau
\end{equation}
\item[$({e})$]  there exists a constant $c>0$ such that  ${\cal E }({u}(s), {v}(s)) + {\cal H}({v}(s))\leq c$ for  $s \in [0,1]$.
\end{description}

We observe that this model is considerably similar to the regularized Francfort \& Marigo model introduced in \cite{MR1745759}.    More precisely, rather than the separate minimality condition $(c)$ the later model uses the global minimality condition
\begin{description}
\item[$(\tilde{c})$]  for all $(\tilde{u},\tilde{v}) \in H^1(\Omega) \times H^1(\Omega)$ with $\tilde{u}= g(s)$ ,  $\tilde{v}=1$ on $\partial \Omega$  and  $0 \leq \tilde{v} \leq v(s)$ we have
\begin{equation}
{\cal E }(u(s), v(s)) + {\cal H}(v(s)) \leq {\cal E }(\tilde{u}, \tilde{v}) + {\cal H}(\tilde{v}).
\end{equation}
\end{description}

We know workout some consequences of the first minimality condition in (c) above, that is, equation \eqref{c.1}.  From definition of subdiferential  we see that
$$
{\cal E}(u,w) + {\cal H}(w) -{\cal E}(u,v) - {\cal H}(v)= \partial_v \Big({\cal E}(u,v) + {\cal H}(v) \Big) (w-v).
$$
Hence,   the first minimality condition in {$(c)$} implies,
\begin{equation}
\partial_v \Big({\cal E}(u,v) + {\cal H}(v) \Big) \leq 0.
\end{equation}
Computing the subdiferential applied to a function $\phi \leq 0 $ gives, 
\begin{multline}
\frac{d}{d\zeta}\frac{1}{2} \int_{\Omega} ((v+\zeta \phi)^2+ \eta_\gamma) |\nabla u|^2dx + \int_{\Omega} \frac{1}{4\gamma} (1-(v+\zeta \phi))^2 + \gamma |\nabla (v+\zeta \phi)|^2 dx \\=\int_{\Omega} (v+\zeta \phi)\phi  |\nabla u|^2dx
+ \int_{\Omega} \frac{-1}{2\gamma} (1-(v+\zeta \phi))\phi  + \gamma 2 (\zeta |\nabla \phi|^2 +   \nabla v\cdot \nabla \phi)  dx
\end{multline}
and by taking $\zeta=0$ yields
$$
\int_{\Omega} v\phi  |\nabla u|^2dx
+ \int_{\Omega} \frac{-1}{2\gamma} (1-v)\phi  + \gamma 2  \nabla v\cdot \nabla \phi  dx \geq  0.
$$
Integrating by parts the last term we obtain
\begin{equation}\label{star}
-2\gamma \Delta v +v |\nabla u|^2    \leq  
 \frac{1}{2\gamma} (1-v).
\end{equation}


\section{Asymptotic expansions}\label{sec:expansions}

Let  $B$ represent a  closed ball of radius $r$ centered at  $(1/2,..., 1/2)$ such that $B \subset (0,1)^d$ and $\Omega \subset \mathbb{R}^d$ be an open set. Let $Y= (0,1)^d \setminus B$ and $\Gamma= \partial Y \cap \partial B$. Also, define  $\Omega_\epsilon$ as the region obtained by the  periodic extension of the  rescaled domain $\epsilon Y$   over $\mathbb{R}^d$ that is contained by $\Omega$.

\begin{figure}\label{fig:omega}
\includegraphics[scale=0.25]{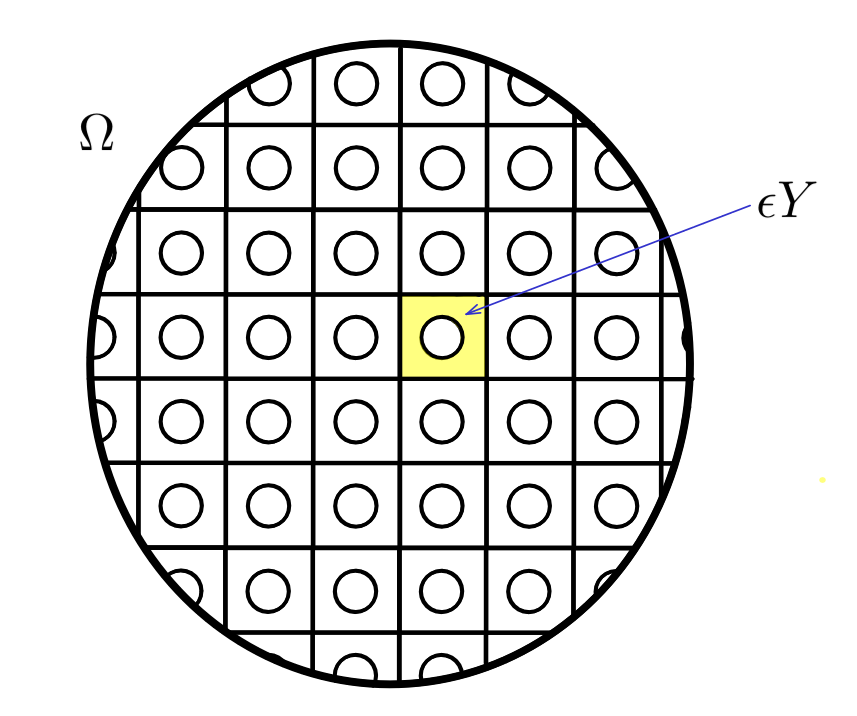}
\caption{Illustration of the domain $\Omega$ and periodic extension of the 
rescaled domain $\e Y$.}
\end{figure}

We set $u_\epsilon$ and $v_\epsilon$ as the solution of the model described in Section \ref{sec:models} applied to the domain $\Omega_\epsilon$ considering that the boundary condition from item $(a)$ is applied in  $ \partial \Omega \cap \partial \Omega_\epsilon$ and we impose the following boundary conditions  on $\Gamma_\epsilon= \partial \Omega_\epsilon \setminus \partial \Omega$, 
\begin{equation}\label{cond-f-G}
\vec{n} \cdot \nabla u_\epsilon= 0 \mbox{~ and ~}  \vec{n}  \cdot \nabla v_\epsilon= 0 \mbox{~ at ~} \Gamma_\epsilon
\end{equation}
where $\vec{n}$ denotes the normal on $\Gamma_\epsilon$ exterior to $\Omega_\epsilon$.
We now use the {\it anzatz}
\begin{eqnarray}\label{ueps}
u_\epsilon &=& u_0(x,y)+\e u_1(x,y)+\e^2 u_2(x,y)+ \cdots \mbox{ and } \\ v_\epsilon
&=& v_0(x,y)+\e v_1(x,y)+\e^2 v_2(x,y)+\cdots \nonumber
\end{eqnarray}
where for each $x \in \Omega $ the functions $u_i(x,\cdot)$ and $v_i(x,\cdot)$ ($i \in\{0,1,2...\}$) are $Y$-periodic. Note that in the following, unless is necessary, we do not make explicit reference to the pseudo-time variable $s$.

\subsection{Minimality condition with respect to $v$}
It is instructive to consider the asymptotic expansion in (\ref{star}) 
to see the consequence of the minimality condition with respect to $v$ in $({c})$.
By replacing the expansions of $u_ \epsilon$ and $v_\epsilon$ in \eqref{star}
we obtain,
\begin{equation}
-2\gamma \Delta v_\e +v_\e |\nabla u_\e|^2    \leq  
 \frac{1}{2\gamma} (1-v_\e)
\end{equation}
or, in terms of asymptotic expansions in the parameter $\e$,  
\begin{multline}\label{des-asympt-expv0}
 -2\gamma (\nabla_x +\frac{1}{\e}\nabla_y) \cdot \left(\nabla_x +\frac{1}{\e}\nabla_y\right) (v_0+ \e v_1+ \e^2v_2+\cdots )  \\ +  (v_0+ \e v_1+ \e^2v_2+\cdots )\left(  
\left| \left(\nabla_x +\frac{1}{\e}\nabla_y\right) (u_0+\e u_1+\e^2 u_2+ \cdots ) \right|^2 \right)\\
\leq  \frac{1}{2\gamma} (1-v_0- \e v_1-\e^2v_2-\cdots ).
\end{multline}
We now work out some terms in this last 
inequality that we need in the rest of the paper. 
Consider first the product  $|\nabla u_\epsilon|^2$. We have, 
\begin{eqnarray}\label{gradu}
|\nabla u_\epsilon|^2&=&\left|\left(\nabla_x +\frac{1}{\e}\nabla_y\right)(u_0+\e u_1+\e^2 u_2+ \cdots )\right|^2\\
&=&\frac{1}{\e^2}{A}_{-2}+\frac{1}{\e}{A}_{-1}+{A}_0+\e {A}_1 + \cdots, \nonumber
\end{eqnarray}
where
\begin{eqnarray*}   
{A}_{-2}&=&|\nabla_y u_0|^2,\\ 
{A}_{-1}&=&2\langle \nabla_y u_0,\nabla_x u_0  \rangle, \\ 
{A}_0&=& |\nabla_x u_0|^2 
+ 2\langle \nabla_y u_1,\nabla_x u_0  \rangle + |\nabla_y u_1|^2 +2\langle \nabla_y u_0,\nabla_x u_1  \rangle + 2\langle \nabla_y u_0,\nabla_y u_2  \rangle, \\
{A}_0&=& |\nabla_x u_0+\nabla_y u_1|^2 +2\langle \nabla_y u_0,\nabla_x u_1 +\nabla_y u_2  \rangle, \\ 
 {A}_1&=&  2\Big[\langle \nabla_x u_1,\nabla_x u_0  \rangle +\langle \nabla_x u_1,\nabla_y u_1  \rangle +\langle \nabla_x u_0,\nabla_y u_2  \rangle+ \langle \nabla_y u_0,\nabla_x u_2  \rangle\Big].
\end{eqnarray*}

We now consider the term $v_\epsilon |\nabla u_\epsilon|^2$.
We have, 
\begin{eqnarray}\label{vgradu}
v_\epsilon |\nabla u_\epsilon|^2&=& (v_0+ \e v_1+ \e^2v_2+\cdots )\left(\frac{1}{\e^2}{A}_{-2}+\frac{1}{\e}{A}_{-1}+{A}_0+\e {A}_1 + \cdots\right)\\\nonumber
&=&\frac{1}{\e^2} v_0{A}_{-2}+ \frac{1}{\e}\Big[ v_1{A}_{-2}+ v_0  {A}_{-1}\Big] 
+ \Big[ v_0 {A}_0 +v_1{A}_{-1}+{ v_2{A}_{-2}} \Big]
\\\nonumber
&&+ \e\Big[v_3 {A}_{-2} +v_2{A}_{-1}+v_1{A}_0  \Big]  +\cdots. 
\end{eqnarray}

Additionally, for the first term on the left hand side of \eqref{des-asympt-expv0} we observe that
\begin{eqnarray}\label{laplace}
\Delta v_\epsilon &=& 
(\nabla_x +\frac{1}{\e}\nabla_y) \cdot \left(\nabla_x +\frac{1}{\e}\nabla_y\right) (v_0+ \e v_1+ \e^2v_2+\cdots )\\\nonumber
&=&\frac{1}{\e^2} \Delta_y v_0 +\frac{1}{\e}\Big[2\nabla_x \cdot \nabla_y v_0  +\Delta_y v_1\Big]\\\nonumber &&
+ \Big[\Delta_x v_0 +\Delta_y v_2
+  { 2 \nabla_x \cdot \nabla_y v_0} {+2\nabla_y\cdot\nabla_x v_1} 
\Big]\\\nonumber && +\e \Big[ \Delta_x v_1+ \Delta_y v_3   \Big]+
\cdots.
\end{eqnarray}

We can now use  (\ref{gradu}), 
(\ref{vgradu}) and  (\ref{laplace}) together with  \eqref{des-asympt-expv0} to get an 
asymptotic expansion inequality,  
\begin{equation}\label{finalInequality}
\frac{1}{\e^2} I_{-2}+\frac{1}{\e}I_{-1}+I_0+\e I_1 + \cdots\leq 0,
\end{equation}
where
\begin{eqnarray*}   
I_{-2}&=&-2\gamma \Delta_yv_0  +v_0{A}_{-2},\\
I_{-1}&=&-2\gamma \Big[2\nabla_x \cdot \nabla_y v_0  +\Delta_y v_1\Big]+\Big[ v_1{A}_{-2}+ v_0  {A}_{-1}\Big],  \\
I_{0}&=& -2\gamma \Big[\Delta_x v_0 +\Delta_y v_2
+ 2 \nabla_x \cdot \nabla_y v_0+2\nabla_y\cdot\nabla_x v_1\Big]\\&&+\Big[ v_0 {A}_0 +v_1{A}_{-1}+ {v_2{A}_{-2}} \Big]-  \frac{1}{2\gamma}(1-v_0),\\
I_{1}&=& -2\gamma \Big[ \Delta_x v_1+ \Delta_y v_3   \Big]+\Big[v_3 {A}_{-2} +v_2{A}_{-1}+v_1{A}_0  \Big]-\frac{1}{2\gamma}v_1.
\end{eqnarray*}

Relation \eqref{des-asympt-expv0} is an inequality and therefore we are not allowed to simply match the terms of the same order with respect to $\e$ in this case. We complement the information in \eqref{des-asympt-expv0} by studying 
the balance of energy relation in {\eqref{d}}.

\subsection{Balance of energy}
In order to bring the balance of energy relation in 
{\eqref{d}} we compute 
the requires terms: $\mathcal{E}$, $\mathcal{H}$ and the integral term. See \eqref{d}.\\ 

To work out the term with  $\mathcal{E}$  in {\eqref{d}} we use \eqref{gradu} and observe that
\begin{eqnarray*}
(v_\e^2+\eta_\gamma)|\nabla u_\epsilon|^2&=&
(v_0^2+\eta_\gamma +\e 2v_0v_1 + 2\e^2v_0v_2 + \e^2v_1^2 + \cdots )\\&& \left(\frac{1}{\e^2}{A}_{-2}+\frac{1}{\e}{A}_{-1}+{A}_0+\e {A}_1 + \cdots\right)\\
&=&\frac{1}{\e^2}(v_0^2+{  \eta_\gamma}){A}_{-2}+\frac{1}{\e}\Big[(v_0^2+{\eta_\gamma} ){A}_{-1}+2v_0v_1{A}_{-2} \Big]\\
&&+\Big[(v_0^2+\eta_\gamma){A}_0+ { 
2v_0v_1{A}_{-1}+2v_0v_2{A}_{-2}} \Big]\\&& +
\e \Big[2v_0v_1 {A}_0+{  2v_0v_2{A}_{-1} +\cdots }  \Big]+\cdots\\
&=&\frac{1}{\e^2}{C}_{-2}+\frac{1}{\e}{C}_{-1}+{C}_0+\e {C}_1 + \cdots
\end{eqnarray*}
where we have introduced the notation ${C}_i$, $i=-2,1,0,1, \dots$ given by 
\begin{eqnarray*}   
{C}_{-2}&=&(v_0^2+{  \eta_\gamma}){A}_{-2}=(v_0^2+{  \eta_\gamma})|\nabla_y u_0|^2,\\ 
{C}_{-1}&=&(v_0^2+{  \eta_\gamma} ){A}_{-1}+2v_0v_1{A}_{-2}, \\ 
{C}_0&=& (v_0^2+\eta_\gamma){A}_0+ { 
2v_0v_1{A}_{-1}+2v_0v_2{A}_{-2}} , \\ 
 {C}_1&=& 2v_0v_1 {A}_0+{  2v_0v_2{A}_{-1} +\cdots }.
\end{eqnarray*}
This last expansion allow us to write the energy $\mathcal{E}$ as, 
\begin{eqnarray} \label{energyE}
{\cal E}(u_\epsilon,v_\epsilon)&=& \frac{1}{2} \int_{\Omega_\epsilon} (v_\e^2+\eta_\gamma)|\nabla u_\epsilon|^2dx 
\\ \nonumber 
&=& { \frac{1}{2}} \int_{\Omega_\epsilon}\left(\frac{1}{\e^2}{C}_{-2}+\frac{1}{\e}{C}_{-1}+{C}_0+\e {C}_1 + \cdots \right) dx. 
\end{eqnarray}

Now we compute $\mathcal{H}$ in  \eqref{d}. As before we write 
\begin{eqnarray}\label{gradv}
|\nabla v_\epsilon|^2&=&\left|\left(\nabla_x +\frac{1}{\e}\nabla_y\right)(u_0+\e u_1+\e^2 u_2+ \cdots )\right|^2\\
&=&\frac{1}{\e^2}{B}_{-2}+\frac{1}{\e}{B}_{-1}+{B}_0+\e {B}_1 + \cdots, \nonumber
\end{eqnarray}
where
\begin{eqnarray*}   
 {B}_{-2}&=&|\nabla_y v_0|^2,\\ 
 {B}_{-1}&=&2\langle \nabla_y v_0,\nabla_x v_0  \rangle, \\ 
 {B}_0&=& |\nabla_x v_0|^2 
+ 2\langle \nabla_y v_1,\nabla_x v_0  \rangle + |\nabla_y v_1|^2 +2\langle \nabla_y v_0,\nabla_x v_1  \rangle + 2\langle \nabla_y v_0,\nabla_y v_2  \rangle, \\ 
{B}_0&=& |\nabla_x v_0+\nabla_y v_1|^2 +2\langle \nabla_y v_0,\nabla_x v_1 +\nabla_y v_2  \rangle, \\
 {B}_1&=&  2\Big[\langle \nabla_x v_1,\nabla_x v_0  \rangle +\langle \nabla_x v_1,\nabla_y v_1  \rangle +\langle \nabla_x v_0,\nabla_y v_2  \rangle+ \langle \nabla_y v_0,\nabla_x v_2  \rangle\Big].
\end{eqnarray*}
We also can write, 
\begin{eqnarray*}
(1-v_\e)^2&=&(1-v_0 -\e v_1- \e^2v_2-\cdots )^2\\&=&
 (1-v_0)^2  -\e 2v_1(1-v_0) + \e^2v_1 -2\e^2 v_2(1-v_0) +\cdots .
\end{eqnarray*}
With these last two expansions we can write, 
 \begin{eqnarray}\label{energyH}
&&{\cal H}(v_\epsilon)=\int_{\Omega_\epsilon} \frac{1}{4\gamma} (1-v_\e)^2+\gamma | \nabla v_\e|^2 dx \\\nonumber
&&=\int_{\Omega_\epsilon} 
\frac{1}{\e^2} \gamma {B}_{-2} + 
\frac{1}{\e} \gamma {B}_{-1}+\Big[
\frac{1}{4\gamma} (1-v_0)^2  +\gamma {B}_0 \Big]+\\\nonumber
&&+
\e \Big[-\frac{1}{4\gamma}2v_1(1-v_0) +\gamma{B}_1
\Big]+\e^2\Big[- \frac{1}{4\gamma}(v_2(1-v_0) +\cdots)+\gamma{B}_2+\cdots  \Big]dx. 
\end{eqnarray}

Finally, in order to use {\eqref{d}} we compute the term 
$\displaystyle\int_0^s \langle  (\eta_\gamma + {v}^2_\e) \nabla {u}_\e, \nabla g_s \rangle d\tau$. For this we write, 
\begin{eqnarray}\label{auxtermvepssquaregradu}
 &&(v^2_\epsilon+\eta_\gamma) \nabla u_\epsilon =
 (v_0^2+\eta_\gamma +\e 2v_0v_1 + 2\e^2v_0v_2 + \e^2v_1^2 + \cdots )\\\nonumber &&\left(\nabla_x +\frac{1}{\e}\nabla_y\right)(u_0+\e u_1+\e^2 u_2+ \cdots )\\\nonumber 
 &&=\frac{1}{\epsilon}  (v_0^2 + \gamma) \nabla_y u_0  + (v_0^2 + \gamma) (  \nabla_x u_0 + \nabla_y u_1)+2v_0 v_1\nabla_y u_0  \\\nonumber  &&+
\e \Big[  2v_0v_1  (\nabla_x u_0 + \nabla_y u_1)+ (v_0^2 + \gamma)( \nabla_x u_1 + \nabla_y u_2)  \\\nonumber &&+  (2v_0  v_2+ v_1^2) \nabla_y u_0 \Big] + \cdots \\\nonumber
&&=\frac{1}{\e} {D}_{-1}+{D}_0+\e {D}_1+\cdots
\end{eqnarray}
where we have introduced the notation ${D}_i$, $i=-1,0,1,\dots$ given by
\begin{eqnarray*}   
{D}_{-1}&=&  (v_0^2 + \gamma) \nabla_y u_0, \\
{D}_{0}&=& (v_0^2 + \gamma) (  \nabla_x u_0 + \nabla_y u_1)+2v_0 v_1\nabla_y u_0,  \\
{D}_{1}&=&  2v_0v_1  (\nabla_x u_0 + \nabla_y u_1)+ (v_0^2 + \gamma) (\nabla_x u_1 + \nabla_y u_2) + \\&& (2v_0  v_2+ v_1^2)\nabla_y u_0.
\end{eqnarray*}
Therefore, we can write 
\begin{multline}\label{integralterm}
\int_0^s \Big\langle  (\eta_\gamma + v_\e^2) \nabla u_\e, \nabla g_\tau \Big\rangle d\tau= \int_0^s \Big\langle\frac{1}{\e} {D}_{-1}+{D}_0+\e {D}_1+\cdots,\nabla g_\tau \Big\rangle    d\tau. \\
=\int_0^s  \frac{1}{\e}\langle {D}_{-1},\nabla g_\tau \rangle +\langle {D}_0,\nabla g_\tau \rangle +\e \langle{D}_1,\nabla g_\tau \rangle  +\cdots d\tau.
\end{multline}

In the conservation of energy relation {\eqref{d}} we insert the expansion of the energy $\mathcal{E}$ in 
\eqref{energyE}, the energy $\mathcal{H}$ in 
\eqref{energyH} and the integral term in 
\eqref{integralterm}. We obtain, 
\begin{equation}\label{finalEnergyConservation}
\frac{1}{\e^2} E_{-2}+\frac{1}{\e}E_{-1}+E_0+\e E_1 + \cdots=0,
\end{equation}
where
\begin{eqnarray*}   
E_{-2}&=& \int_{\Omega_\epsilon} \frac{1}{2}{C}_{-2}(s)+\gamma {B}_{-2}(s)- {C}_{-2}(0)-\gamma {B}_{-2}(0), \\
E_{-1}&=&   \int_{\Omega_\epsilon} \frac{1}{2}{C}_{-1}(s)+\gamma {B}_{-1}(s)-{C}_{-1}(0)-\gamma {B}_{-1}(0)- \int_0^s\langle {D}_{-1},\nabla g_\tau\rangle,  \\
E_{0}&=&  \int_{\Omega_\epsilon} \frac{1}{2}{C}_0(s) +\Big[ \frac{1}{4\gamma}( 1-v_0(s)  )^2+\gamma {B}_0(s)\Big], \\
&&-\frac{1}{2}{C}_0(0)-\Big[ \frac{1}{4\gamma} (1-v_0(0))^2  +\gamma {B}_0(0)\Big]-\int_0^s  \langle {D}_0,\nabla g_\tau \rangle.
\end{eqnarray*}
 
\subsection{Minimality condition with respect to $u$ }
We still need one more ingredient before 
starting the asymptotic expansion comparisons. 
The second minimality condition in $({c})$, equation 
\eqref{c.2}, yields the following equation
\begin{equation}\label{eq-asympt-expu0}
\nabla \cdot [ (v_\epsilon^2+ \eta_\gamma) \nabla u_\epsilon] =0.
\end{equation}
Recalling the expansion for $(v_\e^2+\eta_\gamma)\nabla u_\e$ in 
\eqref{auxtermvepssquaregradu} we get 

$$
 (\nabla_x+ \frac{1}{\e} \nabla_y)  \cdot \Big[ \frac{1}{\e} {D}_{-1}+{D}_0+\e {D}_1+\cdots\Big] =0
$$
and therefore
\begin{equation}\label{FinalMinimality}
    \frac{1}{\e^2}\Big[\nabla_y\cdot{D}_{-1} \Big] + 
    \frac{1}{\e}\Big[ \nabla_x\cdot {D}_{-1}+\nabla_y\cdot{D}_0\Big] 
    + \Big[\nabla_x\cdot {D}_0+\nabla_y\cdot {D}_1  \Big]  + \cdots =0.
\end{equation}

\subsection{Asymptotic expansion terms matching}
We now match the terms with same order with respect to $\e$. 
Recall that the use of asymptotic expansions for $u_\e$ and $v_\e$ led us to:
\begin{itemize}
  \item Balance of energy: equation \eqref{finalEnergyConservation} obtained from the
 balance of energy relation {\eqref{d}},
 \item Minimality condition with respect to $v$: inequality \eqref{finalInequality} obtained from the inequality \eqref{star}, 
 \item Minimality condition with respect to $u$: equation \eqref{FinalMinimality} obtained from \eqref{c.2}.

\end{itemize}

\subsubsection{Terms corresponding to $\e^{-2}$} 
We  consider the terms of order $\e^{-2}$.\\

\underline{\it Balance of energy:} From the term  with order $\e^{-2}$ in equation \eqref{FinalMinimality} and the definition of ${D}_{-1}$ we have
$$
 \nabla_y \cdot  (v_0^2 + \eta_\gamma)   \nabla_y u_0=0.
$$
Since  $\eta_\gamma>0$ and $v_0$ and  $u_0$ is $Y$-periodic with respect to $y$  we conclude that $u_0$ is independent of $y$, that is,
\begin{equation}\label{grad_u0_0}
\nabla_y u_0=0.
\end{equation}
This implies that 
\begin{itemize}
    \item ${A}_{-2}=|\nabla_y u_0|^2=0 $,
    \item ${A}_{-1}=2\langle \nabla_y u_0,\nabla_x u_0  \rangle=0$,
    \item ${D}_{-1}=  (v_0^2 + \gamma) \nabla_y u_0 =0$,
    \item  ${C}_{-2}=(v_0^2+ \eta_\gamma){A}_{-2}=0 $,
    \item ${C}_{-1}=(v_0^2+{  \eta_\gamma} ){A}_{-1}+2v_0v_1{A}_{-2}$=0. 
\end{itemize}

Using ${C}_{-2}=0$ and considering the term of order $\e^{-2}$ in  the balance of energy relation  \eqref{finalEnergyConservation} and \eqref{gradv} we conclude that for all $s\in[0,1]$,
\begin{equation}
\int_{\Omega_\epsilon} |\nabla_y v_0 (s)|^2 \;dx = \int_{\Omega_\epsilon}  |\nabla_y v_0 (0)|^2\;dx.
\end{equation}
From condition $({a})$  and $(b)$ of 
the regularized quasi-static fracture model $v_\e$ is a decreasing function with respect to the $s$ variable, also $0\leq v_\epsilon\leq 1$.
By taking limits when $\e\to 0$ we see that these properties also are valid for 
$v_0$: it is a decreasing function with respect to the $s$ variable and  $0\leq v_0\leq 1$.
We expect
$v_\epsilon$ to be close to one in most of the domain $\Omega_\epsilon$.  Using this as our motivation we choose  $v_0$ independent of $y$, that is, 
\begin{equation}\label{grad_v0_0}
\nabla_y v_0=0.
\end{equation}
This implies that 
\begin{itemize}
    \item ${B}_{-2}=|\nabla_y v_0|^2=0 $,
    \item ${B}_{-1}=2\langle \nabla_y v_0,\nabla_x v_0  \rangle =0$. \\
\end{itemize}

\underline{\it Minimality condition with respect to $v$:}
 Multiplying inequality \eqref{finalInequality} by $\e^{2}$, taking $\epsilon \rightarrow 0$ and using
  (\ref{gradu}), 
(\ref{vgradu}) , (\ref{laplace}),  
   the last equation and \eqref{grad_v0_0},  we obtain
   $$
I_{-2}=-2\gamma \Delta_yv_0  +v_0{A}_{-2}\leq 0.
$$
If we recall that  ${A}_{-2}=|\nabla_y u_0|^2=0$ we get
$$
- \Delta_y v_0 \leq 0.
$$
We observe that this relation is satisfied since we  already obtained that $v_0$ is independent of $y$; see \eqref{grad_v0_0}.

\subsubsection{Terms corresponding to $\e^{-1}$}
We now consider the terms of order $\e{-1}$.\\

\underline{\it Balance of energy:}
The term of order $\e^{-1}$ obtained from balance of energy relation is (see \eqref{finalEnergyConservation})
$$
E_{-1}=  \int_{\Omega_\epsilon}  \frac{1}{2}{C}_{-1}(s)+\gamma {B}_{-1}(s)-{C}_{-1}(0)-\gamma {B}_{-1}(0)-  \int_0^s\langle {D}_{-1},\nabla g_\tau\rangle =0
$$
which is satisfied due to 
${C}_{-1}=0$, ${B}_{-1}=0$ and ${D}_{-1}=0$ as we have seen before.\\

\underline{\it Minimality condition with respect to $v$:} Since  $\nabla_y v_0=0$ and $\nabla_y u_0=0$ the   term of order $\epsilon^{-2}$ in inequality \eqref{finalInequality}, $I_{-2}$, vanishes. Hence,   we can multiply \eqref{finalInequality} by $\e>0$ and take $\epsilon \to 0^+$ to obtain 
$$
-2\gamma \Delta_y v_1  \leq 0 ~\mbox{ for  }  y\in Y.
$$


We observe  that if $-\Delta_y v_1 < 0$  for some $x \in \Omega$ by taking the limit $\epsilon \rightarrow 0$ in the asymptotic expansion \eqref{des-asympt-expv0} we are not able  to conclude any information about the  zero order term with respect to $\epsilon$.  Since we want to obtain a model based on the limit of  asymptotic expansion  we impose the more restrictive  condition
\begin{equation}\label{delta_v1_0}
-\Delta_y v_1  =0 ~\mbox{ for  }  y\in Y.
\end{equation}
\begin{remark}
As an extra motivation for this choice note that if we where allow to use $\e<0$, then multiplying 
\eqref{finalInequality} by $\e$ and taking $\e\to 0^-$ would give us 
$-2\gamma \Delta_y v_1  \geq 0$.\\
\end{remark}
\underline{\it Minimality condition with respect to $u$:} matching the terms of  order  $\e^{-1}$   in equation \eqref{FinalMinimality} we obtain 
$$
\nabla_x\cdot {D}_{-1}+\nabla_y\cdot{D}_0=0
$$ that after replacing ${D}_{-1}=0$ and 
${D}_0$ give as
\begin{equation}\label{term_eps-1_exp_u}
 \nabla_y\cdot \Big[ (v_0^2 + \gamma) (  \nabla_x u_0 + \nabla_y u_1)
 \Big]=0.
\end{equation}
Since $u_0$ and $v_0 $ are independent of $y$, see \eqref{grad_u0_0} and \eqref{grad_v0_0},  equation \eqref{term_eps-1_exp_u} can be rewritten as
\begin{equation}\label{delta_u1_0}
\Delta_y u_1=0 ~\mbox{ for } y \in Y.
\end{equation}

\subsubsection{Terms corresponding to $\e^{0}$} This time we consider terms of order zero. \\

\underline{\it Balance of energy:}
We start with the expansion coming from the balance of energy relations, that is, equation \eqref{finalEnergyConservation}. From the term  of order $\e^0$ in the conservation of energy relation we obtain 
$$
E_0=0
$$
or 
\begin{eqnarray*}
&&{ } \int_{\Omega_\epsilon}  \frac{1}{2}{C}_0(s) +\Big[ \frac{1}{4\gamma}( 1-v_0(s)  )^2+\gamma {B}_0(s)\Big]
\\&&\quad- \frac{1}{2}{C}_0(0)-\Big[ \frac{1}{4\gamma} (1-v_0(0))^2   +\gamma {B}_0(0)\Big]-\int_0^s  \langle {D}_0,\nabla g_\tau \rangle=0
\end{eqnarray*}

Using that ${A}_{-2}=0$, ${A}_{-1}=0$, $\nabla_y u_0=0$ and  $\nabla_y v_0=0$ we have
$$
C_0= (v_0^2+\eta_\gamma)|\nabla_x u_0+\nabla_y u_1|^2 , \quad
B_0= |\nabla_x v_0+\nabla_y v_1|^2
$$
and 
$$
D_0=(v_0^2 + \gamma) (  \nabla_x u_0 + \nabla_y u_1).
$$
Therefore and recalling that $\langle f,g\rangle =\int_{\Omega_\epsilon} fg$,  we can write
\begin{equation}\label{order0minilaity}
R_\e(s)=
 \int_0^s \int_{\Omega_\e} (v_0^2(t) + \gamma) (  \nabla_x u_0(t) + \nabla_y u_1(t),\nabla g_\tau(t) \; dx dt
+R_\e(0)
\end{equation}
where 
\begin{eqnarray} 
R_\e(t)&=&
{\int_{\Omega_\e}}\frac{1}{2}(v_0(t)^2+\eta_\gamma)\Big(|\nabla_x u_0(t)+\nabla_y u_1(t)|^2 \Big)\\
&&+ \frac{1}{4\gamma} (1-v_0(t))^2   +\gamma 
\Big(|\nabla_x v_0(t)+\nabla_y v_1(t)|^2\Big)\; dx.
\end{eqnarray}


\underline{\it Minimality condition with respect to $v$:}
Since the term of order $\epsilon^{-2}$  and 
$\e^{-1}$, $I_{-2}$ and $I_{-1}$, are zero in inequality  \eqref{finalInequality}, we can take the limit $\epsilon \rightarrow 0$  to conclude from the zero order term that 
$$
I_0\leq 0.
$$
Recalling that 
\begin{eqnarray*}
I_{0}&=& -2\gamma \Big[\Delta_x v_0 +\Delta_y v_2
+ 2 \nabla_x \cdot \nabla_yv_1 \Big]\\&&\quad +\Big[ v_0 {A}_0 +v_1{A}_{-1}+ {v_2{A}_{-2}} \Big] -\frac{1}{2\gamma}(1-v_0),
\end{eqnarray*}
and that ${A}_{-2}={A}_{-1}=0$  and 
$\nabla_y v_0=0$ we obtain,
$$
-2\gamma \Big[\Delta_x v_0 +\Delta_y v_2+2 \nabla_y\cdot\nabla_x v_1\Big]+\Big[ v_0 {A}_0 \Big] -\frac{1}{2\gamma}(1-v_0)\leq 0.
$$
By substituting ${A}_0$
\begin{eqnarray*}
&&-2\gamma \Big[\Delta_x v_0 +\Delta_y v_2+ 2 \nabla_x \cdot \nabla_yv_1\Big]\\
&&\quad+v_0\Big[
|\nabla_x u_0|^2 
+ 2\langle \nabla_y u_1,\nabla_x u_0  \rangle + |\nabla_y u_1|^2   \Big] -\frac{1}{2\gamma}(1-v_0)\leq 0
\end{eqnarray*}
that we can write as
\begin{equation}\label{des-v0ordem0}
-2\gamma \Big[\Delta_x v_0 +\Delta_y v_2+ 2 \nabla_x \cdot \nabla_yv_1\Big]+v_0\Big[ |\nabla_x u_0+\nabla_y u_1|^2\Big] \leq \frac{1}{2\gamma}(1-v_0).
\end{equation}

\underline{\it Minimality condition with respect to $u$: }
Matching the terms of  order  $\e^{0}$   in equation \eqref{FinalMinimality} we obtain
$$
\nabla_x\cdot {D}_0+\nabla_y\cdot {D}_1=0
$$
that give us,
\begin{eqnarray}\label{term_eps0_exp_u}
&&\nabla_x\cdot \Big[ (v_0^2 + \gamma) (  \nabla_x u_0 + \nabla_y u_1) \Big]+\\\nonumber &&\quad\quad\nabla_y\cdot\Big[
 2v_0v_1  (\nabla_x u_0 + \nabla_y u_1)+ (v_0^2 + \gamma) (\nabla_x u_1 + \nabla_y u_2)
\Big]=0~\mbox{ in }~  Y.
\end{eqnarray}

\subsubsection{Boundary condition} 
With respect to the boundary conditions at $\Gamma_\epsilon$ in  \eqref{cond-f-G} we conclude that
$$
\left[ \nabla_x + \frac{1}{\epsilon} (\nabla_y  u_0+ \epsilon u_1 +  \cdot \cdot \cdot ) \right] \cdot \vec{n} = 0.
$$
From the terms of order $\epsilon^0$ and $\epsilon^1$ we obtain
\begin{equation}\label{cond_frontueps}
 \nabla_x u_0 \cdot \vec{n}  + \nabla_y u_1 \cdot \vec{n} =0 ~\mbox{ and }~  \nabla_x u_1 \cdot \vec{n}  + \nabla_y u_2 \cdot \vec{n} =0   ~\mbox{ on }~ \Gamma,
\end{equation}
 Analogously for $v_\e$ to obtain, 
\begin{equation}\label{cond_frontveps}
 \nabla_x v_0 \cdot \vec{n}  + \nabla_y v_1 \cdot \vec{n} =0 ~\mbox{ and }~  \nabla_x v_1 \cdot \vec{n}  + \nabla_y v_2 \cdot \vec{n} =0   ~\mbox{ on }~ \Gamma.
\end{equation}

\section{Homogenization analysis}\label{sec:homogenization}
We now proceed following the standard asymptotic expansion  calculations for the homogenization of a second order elliptic equation to handle  equations \eqref{delta_u1_0}-\eqref{term_eps0_exp_u}. See \cite{MR1434316} and related references.  We start with problems \eqref{delta_v1_0}-\eqref{des-v0ordem0}.\\

\subsection{Preliminary observations and definitions}
Denote by $e_i$, $i=1,\dots,n$ the 
canonical vector of $\mathbb{R}^n$. For 
$i=1,\dots,n$ we introduce  the cell problem
\begin{equation}\label{cel-prob}
-\Delta_y z_i = 0, ~\mbox{  }   \nabla_y z_i \cdot \vec{n} =e_i\cdot \vec{n}   ~\mbox{ on }~ \Gamma  ~\mbox{ and  }~ z_i \mbox{ is  }~ Y\mbox{-periodic}
\end{equation}
and we write
\begin{equation}\label{def_v1}
v_1(x,y)=-\sum_{i=1}^n z_i \partial_{x_i} v_0 + \tilde{v}_1(x)
\end{equation}
and 
\begin{equation}\label{def_u1}
u_1(x,y)=-\sum_i z_i \partial_{x_i} u_0 + \tilde{u}_1(x).
\end{equation}
where we have introduced the correctors $\tilde{v}_1$ and $\tilde{u}_2$.\\

In order to average   \eqref{des-v0ordem0} along the 
$y$ direction we need some preliminary observations and definitions. 
We use the $Y$-periodicity  of $v_1$, $u_1$,  and $v_2$ to obtain 
\begin{equation}\label{yav1}
\int_Y \Delta_y v_2= \int_{\partial Y} \nabla_y v_2 \cdot  \vec{n} \; d \sigma +  \int_\Gamma \nabla_y v_2 \cdot  \vec{n} \; d \sigma =\int_\Gamma \nabla_y v_2 \cdot  \vec{n} \; d \sigma.
\end{equation}
Note that from \eqref{yav1}  and 
\eqref{cond_frontveps} we see that, 
\begin{equation}\label{cancellation}
\int_Y \Delta_y v_2+	\int_Y \nabla_y \cdot (\nabla_x v_1)=0.
\end{equation}
In order to handle  the term $ |\nabla_y u_1|^2$ we  observe that 
\begin{multline}\label{gradyu1square}
\int_Y | \nabla_y u_1|^2 \; dy= \int_Y \left|\sum_{i=1}^n  \nabla_y  z_i \partial_{x_i} u_0 \right|^2 \; dy= 
\int_Y |J_y Z^t  \nabla_x u_0 |^2 \; dy \\ =   \nabla_x u_0^t \int_Y  J_y Z  J_y Z^t \; dy \; \nabla_x u_0.
\end{multline}
Here $J_y Z$ represents the Jacobian with respect to the $y$ variable of the function $Z(y)=(z_1(y), ..., z_n(y))^t$.\\

We have the following matrix identity 
\begin{equation}\label{Mat_Id}
    \int_Y J_y Z^t  \;dy  =    \int_Y J_y Z J_y Z^t  \;dy
\end{equation}
 Indeed, the weak formulation of problem \eqref{cel-prob} is 
\begin{equation}\label{cel-prob-wf}
\int_Y \nabla_y(z_i -y_i) \cdot \nabla_y \phi \; dy = 0, ~~ \mbox{ for all }\; \phi \in H^1_{per}(Y) .
\end{equation}
By taking $\phi=z_j$ we have that 
$$
\int_Y \nabla_y z_i  \cdot \nabla_y z_j \; dy = \int_Y e_i \cdot \nabla_y z_j   \; dy $$
yielding  \eqref{Mat_Id}.

We also observe from the  definition of $u_1$   \eqref{def_u1} and identity \eqref{Mat_Id}  that 
\begin{multline}\label{nxu0+nyu1}
\int_Y|\nabla_x u_0+\nabla_y u_1|^2 \;dy =   \int_Y|\nabla_x u_0|^2   -2\nabla_x u_0^t     J_y Z^t   \nabla_x u_0  \\    +    
\nabla_x u_0^t   J_y Z J_y Z^t\nabla_x u_0      \;dy =
|Y|  |\nabla_x u_0|^2 -  \nabla_x u_0^t    \int_Y J_y Z^t  \;dy  \nabla_x u_0= \\
  \nabla_x u_0^t  \left(|Y| I  -     \int_Y J_y Z^t  \;dy \right)   \nabla_x u_0
\end{multline}
where $I$ is the identity matrix.
The following  matrix will play an important role in our homogenized model 
\begin{equation}\label{M0}
M_0=     I - \frac{1}{|Y|} \int_Y J_y Z^t  \;dy   = 
  I - \frac{1}{|Y|} \int_Y J_y Z J_y Z^t  \;dy   
\end{equation}
Using the definition of $M_0$ we can rewrite  identity \eqref{nxu0+nyu1} as
\begin{equation}\label{normau_0u_12}
    \int_Y|\nabla_x u_0+\nabla_y u_1|^2 \;dy = |Y| \nabla_x u_0^t  M_0   \nabla_x u_0.
\end{equation}
Proceeding in a similar manner one can obtain that
\begin{equation}\label{normau_0v_12}
    \int_Y|\nabla_x v_0+\nabla_y v_1|^2 \;dy = |Y| \nabla_x v_0^t  M_0   \nabla_x v_0.
\end{equation}
\\

\subsection{Averaging the minimality condition with respect to $v$  }
We now come back to averaging \eqref{des-v0ordem0}. 
Integrating in both sides of inequality \eqref{des-v0ordem0} with respect to $y$ over $Y$ 
\begin{multline}\label{des-v0ordem0-av}
-2\gamma  \int_Y\Big[    \Delta_x v_0 +\Delta_y v_2+ 2 \nabla_x \cdot \nabla_yv_1\Big]+v_0\Big[ |\nabla_x u_0+\nabla_y u_1|^2\Big]  \;dy   \\
\leq  \int_Y \frac{1}{2\gamma}(1-v_0) \; dy.
\end{multline}

Next, we  estimate the sum of the second and the third term on the left-hand-side of the last inequality using \eqref{cancellation}
\begin{equation}
    \int_Y    \Delta_x v_0 +\Delta_y v_2+ 2 \nabla_x \cdot \nabla_yv_1  \; dy =     \int_Y    \nabla_x \cdot \nabla_yv_1  \; dy.
\end{equation}
 Using the definition of $v_1$  \eqref{def_v1},   and \eqref{normau_0u_12}  we obtain

\begin{eqnarray*}
&& -2\gamma\Big[  \Delta_x v_0|Y|  -  \int_Y \nabla_x \cdot (  J_y Z^t \nabla_{x} v_0) \;dy \Big]   +
v_0|Y|\Big[   \nabla_x u_0^t  M_0   \nabla_x u_0 \Big]
\leq \frac{|Y|}{2\gamma}(1-v_0) .
\end{eqnarray*}
Using the definition of the matrix $M_0$ the fact that $v_0$ is independent of $Y$ and dividing the last equation  by $|Y|$ we obtain  
\begin{equation}\label{eq:des_u0v0}
-2\gamma   \nabla_x \cdot (M_0 \nabla_x v_0)    + 
v_0\Big[    M_0   \nabla_x u_0 \cdot  \nabla_x u_0 \Big] \leq \frac{1}{2\gamma}(1-v_0).
\end{equation}

Problems \eqref{order0minilaity} and \eqref{term_eps0_exp_u} are treated in a similar way. 

\subsection{Averaging the minimality condition with respect to $u$  }

Integrating both sides of \eqref{term_eps0_exp_u} with respect to $y$  over $Y$ we obtain
\begin{eqnarray*}
&& \frac{1}{|Y|} \int_Y \nabla_x\cdot \Big[ (v_0^2 + \gamma) (  \nabla_x u_0 + \nabla_y u_1) \Big]\\\nonumber &&+\frac{1}{|Y|} \int_Y \nabla_y\cdot\Big[
 2v_0v_1  (\nabla_x u_0 + \nabla_y u_1)\Big] +
 \frac{1}{|Y|} \int_Y \nabla_y\cdot \Big[(v_0^2 + \gamma) (\nabla_x u_1 + \nabla_y u_2)
\Big]=0.
\end{eqnarray*}
Next we show that the   second  and third term on the left hand side of the last equation  vanish due to the $Y$-periodicity 
of the functions $u_1$ and $u_2$ and relation \eqref{cond_frontueps}.
In fact, since $v_0$ is independent of $y$ 
\begin{eqnarray*}
 \int_Y \nabla_y\cdot \Big[(v_0^2 + \gamma) (\nabla_x u_1 + \nabla_y u_2)\Big]&=& 
(v_0^2 + \gamma) \int_Y \nabla_y \cdot (\nabla_x u_1 + \nabla_y u_2)\\ 
&=&\int_{\partial Y}  (  \nabla_y u_2+ \nabla_x u_1)\cdot \vec{n}=0
\end{eqnarray*}
due to \eqref{cond_frontueps}.  Similarly, we also have
\begin{multline}
    \int_Y \nabla_y\cdot\Big[
 v_0v_1  (\nabla_x u_0 + \nabla_y u_1)\Big] =  \int_{\partial Y} v_0v_1  (\nabla_x u_0 + \nabla_y u_1)\cdot\vec{n} + \\    \int_{\Gamma} v_0v_1  (\nabla_x u_0 + \nabla_y u_1)\cdot\vec{n} = 0
\end{multline}
  where we used the periodicity to obtain that the first term on the right hand side of the first equation is zero and \eqref{cond_frontveps} to obtain that the second term on the right hand side is zero.  
Therefore, from \eqref{def_u1} and recalling that 
$v_0$ is independent of $y$ we obtain,  
$$
  \nabla_x \cdot  (v_0^2 + \eta_\gamma)   \nabla_x u_0+
 \frac{1}{|Y|} \int_Y  \nabla_x \cdot \left[  (v_0^2 + \eta_\gamma)    \left(-\sum_i \nabla_y z_i \partial_{x_i} u_0\right)\right] \;dy =0.
$$
From the definition of the matrix $M_0$ in \eqref{M0} we have
\begin{equation}\label{def_u_0}
\nabla_x \cdot \Big[ (v_0^2 + \eta_\gamma)  M_0 \nabla_x u_0\Big] =0~\mbox{ in }~ \Omega  ~\mbox{ and }~ u_0=g ~\mbox{ on }~ \partial \Omega.
\end{equation}

\subsection{ Averaging the energy balance relation }

We start observing that  the energy balance relation involves an integral with respect to the the variable $x$ in the original formulation. One may wonder the what is the right meaning of the integral in \eqref{order0minilaity}. More precisely, recalling the definition of $R_\e$ 
\begin{eqnarray*} 
R_\e(t)&=&
{\int_{\Omega_\e}}\frac{1}{2}(v_0(t)^2+\eta_\gamma)\Big(|\nabla_x u_0(t)+\nabla_y u_1(t)|^2 \Big)\\
&&+ \frac{1}{4\gamma} (1-v_0(t))^2   +\gamma 
\Big(|\nabla_x v_0(t)+\nabla_y v_1(t)|^2\Big)\; dx.
\end{eqnarray*}
One could think of the last integral in two different ways: $(i)$ one may interpret the function inside the integral as a function of two independent variables $x$ and $y$ and that we are integrating just with respect to the slow varying $x$ variable, or $(ii)$ may interpret the function inside the integral as a function of two  variables $x$ and $y$ and that the integration is in the  {\it diagonal}  $(x,x/\e)$.  We show that these two interpretations of   equation \eqref{order0minilaity} yield the same limiting relation.

In the case we interpret equation \eqref{order0minilaity} as in $(i)$ described above, the   limiting energy balance relation can be obtained by a  method similar to the one used in previous sections. More precisely,  one can integrate  relation \eqref{order0minilaity}    with respect to $y$ in the $Y$ domain and use relations \eqref{normau_0u_12} and \eqref{normau_0v_12} to obtain the desired limit relation.

Next, we use interpretation $(ii)$ of equation \eqref{order0minilaity} to obtain our limiting model.  Starting from the relation  \eqref{order0minilaity}		and using the definition of the functions $v_1$ \eqref{def_v1} and  $u_1$ \eqref{def_u1}  we have

\begin{multline}\label{order0minilaityPart1}
R_\e=\int_{\Omega_\e } \frac{1}{2} (v_0^2+ \eta)\Big[ -2 \nabla_y z_i \partial_{x_i}u_0 \cdot \nabla_x u_0   + |\nabla_y z_i \partial_{x_i}u_0|^2  + |\nabla_x u_0|^2   \Big] dx \\ +\int_{\Omega_\e } \frac{1}{4\gamma} (1-v_0)^2 + \gamma |\nabla_x v_0|^2+ \gamma|\nabla_y z_i \partial_{x_i}v_0|^2 -\gamma 2\nabla_y z_i \partial_{x_i}v_0\cdot \nabla_x v_0
  \;  dx 
\end{multline}
and 
\begin{multline}\label{order0minilaityPart2}
    \int_0^s \int_{\Omega_\e} (v_0^2(t) + \gamma) (  \nabla_x u_0(t) + \nabla_y u_1(t),\nabla g_\tau(t) \; dx dt=\\ \int_0^s  \langle (v_0^2 + \gamma) (\nabla_x u_0 - \nabla_y z_i \partial_{x_i}u_0) ,\nabla g_s \rangle.
\end{multline}


Here, we reinforce that the functions $\nabla_y z_i$ are evaluated at $x/\e$.  The cell functions are periodic and  assuming they are``smooth'',   in the limit $\e \rightarrow 0$ we expect that
$$
\partial_{y_j} z_i(x/\e) \rightharpoonup  \frac{1}{|Y|}\int_Y \partial_{y_j} z_i(y) \;dy $$ and  $$  \partial_{y_j} z_i(x/\e)\partial_{y_k} z_l(x/\e) \rightharpoonup  \frac{1}{|Y|}\int_Y \partial_{y_j} z_i(y)\partial_{y_k} z_l(y) \;dy
$$
weakly-$*$ in $L^{\infty}$. 
 
We know replace \eqref{order0minilaityPart1} and \eqref{order0minilaityPart2} into 
\eqref{order0minilaity} and take the limit $\epsilon \rightarrow 0$. When taking the limit if we use  relation \eqref{Mat_Id} and  the definition of the matrix \eqref{M0}    we obtain
\begin{multline}\label{balancedenergy000}
\int_\Omega \frac{1}{2}(v_0^2(s)+ \eta)\left[  M_0\nabla_x u_0(s) \cdot  \nabla_x u_0 (s)       \right] dx \\
+\int_{\Omega} \frac{1}{4\gamma} (1-v_0(s))^2 + \gamma M_0 \nabla_x v_0(s) \cdot  \nabla_x v_0 (s)
    dx =
\int_0^s  \langle (v_0^2 + \gamma)  M_0 \nabla_x u_0 ,\nabla g_s \rangle \;d\tau \\
 \int_\Omega \frac{1}{2}(v_0^2(0)+ \eta)\left[  M_0 \nabla_x u_0(0) \cdot  \nabla_x u_0 (0)       \right] dx \\
+\int_{\Omega} \frac{1}{4\gamma} (1-v_0(0))^2 + \gamma M_0 \nabla_x v_0(0) \cdot  \nabla_x v_0 (0) \;
    dx.
\end{multline}
 Hence, the homogenized model is given by \eqref{eq:des_u0v0}, \eqref{def_u_0}, and the last conservation of energy identity.

\section{The novel regularized brittle fracture propagation in porous media model}\label{novel}

We now study the minimality condition associated with equations \eqref{eq:des_u0v0} and \eqref{def_u_0}. We start defining the energy functionals
\begin{eqnarray}\label{eng_func_E_H0}
{\cal E}_0(u_0,v_0)&=& \frac{1}{2} \int_{\Omega} (v^2+ \eta_\gamma) (M_0\nabla u_0 )\cdot \nabla u_0 \; dx ~~\mbox{and} \\{\cal H}_0(v_0)&=& \int_{\Omega} \frac{1}{4\gamma} (1-v_0)^2 + \gamma (M_0\nabla v_0) \cdot \nabla v_0 dx
\end{eqnarray}
where the matrix $M_0$ is defined  by   \eqref{M0}. We observe that equation \eqref{def_u_0} corresponds to the minimality condition
\begin{equation}\label{c.20}
{\cal E}_0({u_0}(s),{v}_0(s)) + {\cal H}_0({v_0}(s))= \inf_{\phi-g(s) \in H^1_0(\Omega)}  {\cal E}_{ 0}(\phi,{v_0}(s)) + {\cal H}({v_0}(s)).
\end{equation}
On the other hand, inequality \eqref{eq:des_u0v0} can be derived from the minimality condition
\begin{equation}\label{c.10}
{\cal E}_0({u}_0(s),{v}_0(s)) + {\cal H}_0({v}(s))= \inf_{0 \leq z \leq {v}_0(s)}  {\cal E}_0({u}(s),z) + {\cal H}_0(z).
\end{equation}
Note that \eqref{c.10} and \eqref{c.20} are analogous to equations \eqref{c.1} and \eqref{c.2}.

 Additionally we note that \eqref{balancedenergy000} can be written as 
 \begin{eqnarray}\label{d0}
{\cal E }_0({u}_0(s), {v}_0(s)) + {\cal H}_0({v}_0(s)) &=& {\cal E }({u}(0), {v}(0)) + {\cal H}_0({v}(0)) 
\\\nonumber && + \int_0^s \langle  (\eta_\gamma + {v}_0^2) M_0 \nabla_x {u}_0, \nabla g_s \rangle d\tau
\end{eqnarray}
which corresponds to \eqref{d}.
 
\section{Conclusions and perspectives}   
   Using  homogenization theory and formal asymptotic expansions applied to a regularized brittle fracture model, we were able 
   to obtain a homogenized regularized fracture model for a homogeneous porous media. 
The obtained  homogenized model  it is presented in Section \ref{novel}.  Our novel homogenized model is similar in nature and structure to the model we started our analysis with: a regularized fracture model for homogeneous porous media. The main difference between the initial model and the homogenized model is the appearance homogenized energies where gradients are weighted by a tensor
that depends only on the microstructure. This is similar to the case of porous media flow (Darcy flow) when derived from Stokes flow using similar technique. 
We believe this is a novel and promising proposal for fracture modeling in porous media. The obtained novel model had led us to some interesting theoretical questions as well as questions related to the numerical analysis and the design of numerical methods in the spirit of multiscale finite element methods among other considerations such as  modeling suitableness  associate to  \eqref{eq:des_u0v0} and \eqref{def_u_0}. All these issues are subject of current research and we hope to presented some interesting results in future works.

\section{Acknowledgments} \label{sec:acknow}

J. Galvis thanks the great hospitality of  Universidade Federal de Minas Gerais visiting summer programme. J. Galvis also thanks partial support from the European Union's Horizon 2020 research and innovation programme under the Marie Sklodowska-Curie grant agreement No 777778 (MATHROCKS).

\bibliographystyle{plain}
\bibliography{refe}

\end{document}